\documentclass[12pt]{amsart}
\usepackage{verbatim}
\usepackage{amssymb, amsmath}
\usepackage{graphicx}
\usepackage{caption}
\usepackage{subcaption}
\usepackage{mathtools}
\usepackage{bbm}
\usepackage{hyperref}
\usepackage{pdfpages}

\usepackage{color}

\begin{document}
\newtheorem{thm}{Theorem}[section]
\newtheorem*{thm*}{Theorem}
\newtheorem{lem}[thm]{Lemma}
\newtheorem{prop}[thm]{Proposition}
\newtheorem{cor}[thm]{Corollary}
\newtheorem*{conj}{Conjecture}
\newtheorem{proj}[thm]{Project}
\newtheorem{question}[thm]{Question}
\newtheorem{rem}{Remark}[section]

\theoremstyle{definition}
\newtheorem{defn}[thm]{Definition}
\newtheorem{remark}[thm]{Remark}
\newtheorem{exercise}{Exercise}
\newtheorem*{exercise*}{Exercise}

\numberwithin{equation}{section}

\newcommand{\rad}{\operatorname{rad}}

\newcommand\SLASH{\char`\\}
\newcommand{\Z}{{\mathbb Z}} 
\newcommand{\Q}{{\mathbb Q}}
\newcommand{\R}{{\mathbb R}}
\newcommand{\C}{{\mathbb C}}
\newcommand{\N}{{\mathbb N}}
\newcommand{\FF}{{\mathbb F}}
\newcommand{\fq}{\mathbb{F}_q}
\newcommand{\rmk}[1]{\footnote{{\bf Comment:} #1}}

\renewcommand{\mod}{\;\operatorname{mod}}
\newcommand{\ord}{\operatorname{ord}}
\newcommand{\TT}{\mathbb{T}}
\renewcommand{\i}{{\mathrm{i}}}
\renewcommand{\d}{{\mathrm{d}}}
\renewcommand{\^}{\widehat}
\newcommand{\HH}{\mathbb H}
\newcommand{\Vol}{\operatorname{vol}}
\newcommand{\area}{\operatorname{area}}
\newcommand{\tr}{\operatorname{tr}}
\newcommand{\norm}{\mathcal N} 
\newcommand{\intinf}{\int_{-\infty}^\infty}
\newcommand{\ave}[1]{\left\langle#1\right\rangle} 
\newcommand{\Var}{\operatorname{Var}}
\newcommand{\Prob}{\operatorname{Prob}}
\newcommand{\sym}{\operatorname{Sym}}
\newcommand{\disc}{\operatorname{disc}}
\newcommand{\CA}{{\mathcal C}_A}
\newcommand{\cond}{\operatorname{cond}} 
\newcommand{\lcm}{\operatorname{lcm}}
\newcommand{\Kl}{\operatorname{Kl}} 
\newcommand{\leg}[2]{\left( \frac{#1}{#2} \right)}  
\newcommand{\Li}{\operatorname{Li}}

\newcommand{\sumstar}{\sideset \and^{*} \to \sum}

\newcommand{\LL}{\mathcal L} 
\newcommand{\sumf}{\sum^\flat}
\newcommand{\Hgev}{\mathcal H_{2g+2,q}}
\newcommand{\USp}{\operatorname{USp}}
\newcommand{\conv}{*}
\newcommand{\dist} {\operatorname{dist}}
\newcommand{\CF}{c_0} 
\newcommand{\kerp}{\mathcal K}

\newcommand{\Cov}{\operatorname{cov}}
\newcommand{\Sym}{\operatorname{Sym}}

\newcommand{\Ht}{\operatorname{Ht}}

\newcommand{\E}{\operatorname{\mathbb E}} 
\newcommand{\sign}{\operatorname{sign}} 
\newcommand{\meas}{\operatorname{meas}} 
\newcommand{\length}{\operatorname{length}} 

\newcommand{\divid}{d} 

\newcommand{\GL}{\operatorname{GL}}
\newcommand{\SL}{\operatorname{SL}}
\newcommand{\re}{\operatorname{Re}}
\newcommand{\im}{\operatorname{Im}}
\newcommand{\res}{\operatorname{Res}}
 \newcommand{\eigen}{\Lambda} 
\newcommand{\tens}{\mathbf t} 
\newcommand{\diam}{\operatorname{diam}}
\newcommand{\fixme}[1]{\footnote{Fixme: #1}}
 \newcommand{\EWp}{\mathbb E^{\rm WP}} 
\newcommand{\orb}{\operatorname{Orb}}
\newcommand{\supp}{\operatorname{Supp}}
\newcommand{\mmfactor }{\textcolor{red}{c_{\rm Mir}}}
\newcommand{\Mg}{\mathcal M_g} 
\newcommand{\MCG}{\operatorname{Mod}} 
\newcommand{\Diff}{\operatorname{Diff}} 
\newcommand{\If}{I_f(L,\tau)}
\newcommand{\SigGOE}{\Sigma^2_{\rm GOE}}

\title[No PPC for sums of two squares]{No Poissonian Pair Correlation for the sum of two squares}

\author{Sharon Zilberhertz }
\address{School of Mathematical Sciences, Tel Aviv University, Tel Aviv 69978, Israel}
\email{Zilberhertz@mail.tau.ac.il}

\date{\today}
\thanks{ Supported by the Israel Science Foundation {grant No.1881/20}}
\maketitle
\begin{abstract}
We prove that the sums of two squares do not have the metric Poissonian pair correlation property.
\end{abstract}
\section{Introduction and Definitions}

For a strictly increasing sequence of naturals $a_n$ and some $\alpha \in [0,1)$ we say that $a_n\alpha$ is uniformly distributed if for any subinterval $J\subset [0,1)$
\[
\frac{1}{N}\#\{ \{a_n\alpha\}\in J\mid n\leq N \}\xrightarrow[N\to \infty]{}|J|
\]
with $\{\cdot\}$ denoting the fractional part of a real number.\newline

Now, we turn to a somewhat stronger definition of uniform distribution modulo 1, Poissonian pair correlation of a sequence.

\subsection{Metric Poissonian Pair Correlation} 
For a strictly increasing sequence of positive reals $a_n$ we say that it has metric Poissonian pair correlation (MPPC for short) if and only if for a.e $\alpha\in\R$ and all $s\ge 0$ we have
\[
\frac{1}{N}\#\{ \|(a_m - a_n)\alpha\| \leq \frac{s}{N} \mid  1\leq m\neq n\leq N\} \xrightarrow[N\to\infty]{}2s.
\]
\newline
Intuitively, we count the number of Diophantine approximations of $\alpha$ by the set of differences of $(a_n)_{n\in\N}$ i.e. the number of times $(a_n - a_m)\alpha$ has a smaller distance than a scaled factor of $1/N$ (by $s>0$) which is the average gap that one expects to get over $N$ points in the unit interval. So the pair correlation function $T$ essentially measures how uniformly the gaps of $(a_n-a_m)\alpha$ are distributed in scaled intervals.\newline

In the last several decades this theme of Poissonian pair correlation has been studied for various of sequences such as $n^k$ for $k = 1$ and $k\ge 2$ integers, then $n^\theta$ with $\theta > 0$ non-integer where the latter has been shown to be true for all $\alpha \neq 0$ and $\theta < \frac{1}{3}$ in \cite{NT}. For various other sequences such as the mentioned above $n^k\alpha$ it turned out to be difficult to prove the property for specific $\alpha$'s. The difficulty of the deterministic problem made mathematicians turn to the probabilistic problem instead, namely the MPPC which we formulated above. For this case there are several known results, among them the sequence of $n^k$ has been shown not to have this property for $k=1$ and has been proven to have this property for integers $k\ge 2$ by Rudnick and Sarnak (see \cite{RS}), later the sequence $n^\theta$ has been proven to have MPPC as well for any $\theta > 1$ non-integer (see \cite{RT}) and later for $0 < \theta < 1$ (see \cite{ELA}) finishing the case of powers of $n$.\newline

From now on we will consider the case of $a_n$ being naturals only, and thus restrict ourselves to $\alpha\in[0,1)$ instead of $\R$. \newline

We define the function $R_N(v)$ on the naturals ($v>0$) to be
\[
R_N(v) = \#\{ a_m-a_n = v\mid 1\leq n < m\leq N \}
\]
and get an alternative definition to pair correlation by considering the function
\[
T(\alpha, N, s) = \frac{1}{N}\sum_{v \in \mathbb{N}} R_N(v)\mathbf{1}_{\{ \|v\alpha\|\leq \frac{s}{N} \}}
\]
and the definition shall then be that $a_n$ has metric Poissonian pair correlation if and only if
\[
T(\alpha, N, s) \xrightarrow[N\to\infty]{}s
\]
for almost every $\alpha$ for any $s>0$. Note: the limit changed from $2s$ to $s$ because we removed the negative duplicated gaps by taking $a_m-a_n$ with $m>n$.
\newline

The technique used frequently to prove that a natural sequence has MPPC is based on the additive energy of the sequence, where we define the additive energy of a sequence to be
\[
E(a_n, N) = \#\{ a_m - a_n = a_k - a_l \mid 1\leq n,m,k,l \leq N \}.
\]
This term satisfies $N^2 \leq E(a_n, N)\leq N^3$ for any given sequence $a_n$ and has connection to the MPPC property through the next theorems

\begin{thm}[Theorem 6 of \cite{TA}]
    Let $a_n$ be a sequence with $E(a_n, N) \ll \frac{N^3}{(\log N)^C}$ for some big enough $C>1$, then $a_n$ is a MPPC.
\end{thm}
\begin{thm}[Theorem 1 of \cite{LWG}]
    Let $a_n$ be a sequence with $E(a_n, N) \gg N^3$ i.e. maximal order of additive energy, then the $a_n$ is not MPPC and in fact it has no PPC for almost all $\alpha$.
\end{thm}
By these theorems we have a clear red zone and green zone for the additive energy such that the sequence will or will not have MPPC, but notice that there is a gap that remains unknown between $\frac{N^3}{(\log N)^C}$ and $N^3$. This leads of course to the natural question regarding how much further can we expand these red and green zones. There has been several proof of existence for sequences $a_n$ with $E = o(N^3)$, beginning with Bourgain's construction of such sequence (see the Appendix of \cite{CA}), and later in \cite{TFB} the authors constructed a sequence for every additive energy of the order
\[
E(a_n, N) \asymp \frac{N^3}{\log N\log\log N...(\log...\log N).}
\]

In \cite{AW} Walker showed that the primes are not MPPC. It is also known that the primes' additive energy is $E \asymp \frac{N^3}{\log N}$. The latter two papers raised the question whether there is some kind of a Khintchine limit on the order of the additive energy, to which accordingly the sequence is MPPC or not, more precisely whether for a sequence $a_n$ with $E(a_n, N) \sim N^3f(N)$ with $f(N)$ decreasing to $0$, is the MPPC property for $a_n$ determined by the convergence or divergence of the sum
\[
\sum_{n=1}^\infty \frac{f(n)}{n}.
\]
This question has been answered in the negative in \cite{CN} where they constructed a sequence with $E \gg \frac{N^3}{(\log N)^{3/4+\delta}}$ which has the MPPC property. \newline

To summarize, the best known red zone (resp. green zone) is $E \gg N^3$ (resp. $E \ll \frac{N^3}{(\log N)^C}$ for $C>1$ big enough), and the red zone cannot be extended below $\frac{N^3}{(\log N)^{3/4}}$. \newline

In this paper we are going to focus on the sequence of sum of two squares in increasing order, and show it has no MPPC. It is not difficult to show that this sequence has additive energy of order $E \asymp \frac{N^3}{\sqrt{\log N}}$.
\begin{thm}[Main Theorem]\label{SOT}
    Let $a_n$ be the sequence of the sum of two squares in a strictly increasing order, then $a_n$ does not have MPPC.
\end{thm}

The proof of this theorem will exploit the somewhat uniform distribution of the sum of two squares along arithmetic progression of small spacing, i.e. if $A(x)$ is the number of sum of two squares up to $x$ then for appropriate $a$ and $q \leq \log x$
\[A_{q,a}(x) \gg \frac{A(x)}{q}\]
with $A_{q,a}(x)$ being the number of sum of two squares up to $x$ which are congruent to $a$ mod $q$.

With some further work we can prove an analogous result when we replace the sum of two squares with any positive definite binary quadratic form, see \cite{SZ}.

\section{Preliminaries}
In preparation we will require two results
\begin{thm}[Theorem 1 of \cite{JDV}]{\label{thm1}}
Let $f(n)$ be some positive arithmetic function with $f(n) = O(n^{-1})$, denote by $\phi(n)$ the Euler's totient function, then for almost every $\alpha$ there are infinitely many $n$'s with
\[\| n\alpha \| \leq f(n)\]
and $\|n\alpha\| = |r-n\alpha|$ for some $(r,n) = 1$, if and only if the following sum diverges
\[
\sum_{n=1}^\infty f(n)\frac{\phi(n)}{n}.
\]
\end{thm}
Theorem \ref{thm1} due to Vaaler (1978) is a weaker version of the Duffin- Schaeffer conjecture, which states a similar result but removes the need for the condition of $O(n^{-1})$, which in our case does not affect us. The conjecture was recently proven by Dimitris Koukoulopoulos and James Maynard (see \cite{DJ}).
\begin{thm}[Theorem 14.7 of \cite{JH}]{\label{thm2}}
    Let $b'(n)$ be the indicator of the odd naturals which are properly represented by sum of two squares, i.e. $b'(n) = 1$ if and only if $n = r^2+t^2$ with $(r,t)=1$, then for $a$ with $(a,q) = 1$ and $a \equiv 1\mod{(4,q)}$, we have
    \[\sum_{\substack{n\leq x\\ n\equiv a(q)}}b'(n) = \frac{c_qx}{q\sqrt{\log x}}(1+O(\frac{\log q}{\log x})^{1/7})\]
    with $c_q \gg 1$. If $q \equiv 0(4)$ then $c_q = 2\kappa\prod_{\substack{p|q\\ p\equiv 3(4)}}(1+\frac{1}{p})$ (with $\kappa = 0.7642\dots$ being Landau-Ramanujan constant). For $(q,4)=1$ it is the same but with factor of $1/4$.
\end{thm}
\section{Lemmas}
\begin{lem}{\label{lem1}}
    Let $B$ be a set with the property that
    \[ \#\{ n\in[M,2M) \mid n\in B\} \gg \frac{M}{(\log\log M)^2} \]
    for any $M\gg 1$. If we fix $h<1$, then for almost all $\alpha$ there exist infinitely many $n\in B$ such that
    \[\|n\alpha\| \leq \frac{1}{n(\log n)^h}.\]
\end{lem}
\subsection{Proof of Lemma \ref{lem1}}
We apply Theorem \ref{thm1} with 
\[f(n) = \begin{cases}

\frac{1}{n(\log n)^h}&\quad n\in B\\
0 &\quad \text{otherwise}

\end{cases}\]
and divide the sum in the theorem into dyadic intervals to get that using $\phi(n) \gg \frac{n}{\log\log n}$ we have
\begin{multline*}
    \sum_{n\in B}f(n)\frac{\phi(n)}{n} \gg \sum_{t=1}^\infty \sum_{\substack{n\in[2^t,2^{t+1}) \\ n\in B}}f(n)\frac{\phi(n)}{n} \gg \sum_{t=1}^\infty \sum_{\substack{n\in[2^t,2^{t+1}) \\ n\in B}}\frac{f(n)}{\log\log n} \\ \gg \sum_{t=1}^\infty \frac{1}{2^{t+1}(\log 2^{t+1})^h\log\log 2^{t+1}}\sum_{\substack{n\in[2^t,2^{t+1}) \\ n\in B}}1,
\end{multline*}
and by the condition on $B$ in the lemma, we have that the inner sum is
\[\gg \frac{2^t}{(\log\log 2^t)^2}.\]
Therefore, we have
\[\sum_{n\in B}f(n)\frac{\phi(n)}{n}\gg \sum_{t=1}^\infty \frac{1}{(\log 2^{t})^h(\log\log 2^{t})^3} \gg  \sum_{t=1}^\infty \frac{1}{t^h(\log{t})^3} = \infty.\]
Now apply Theorem \ref{thm1} to conclude that for almost every $\alpha$ there are infinitely many $n\in B$ with $\|n\alpha\| \leq f(n)$.
\begin{lem}{\label{lem2}}
Let b(n) be the indicator function for the sum of two squares, then for odd $q \ll \log x$ and natural $a$ such that $(a,q) = 1$ we have for $x\ge y \gg x$
\[\sum_{\substack{x\leq n\leq x+y\\ n \equiv a(q)}}b(n) \gg \frac{y}{q\sqrt{\log x}}.\]
\end{lem}
\subsection{Proof of Lemma \ref{lem2}}
Let $a, q, x, y$ be as in the lemma, and denote by $b(n)$ (resp. $b'(n)$) the naturals represented by the sum of two squares (resp. odd and properly represented by the sum of two squares) then by Theorem \ref{thm2} we have
\begin{multline*}
\sum_{\substack{x < n\leq x+y\\ n\equiv a(q)}}b(n) \ge \sum_{\substack{x < n\leq x+y\\ n\equiv a(q)}}b'(n) = \sum_{\substack{n\leq x+y\\ n\equiv a(q)}}b'(n) - \sum_{\substack{n\leq x\\ n\equiv a(q)}}b'(n) \\= \frac{c_q(x+y)}{q\sqrt{\log(x+y)}} - \frac{c_qx}{q\sqrt{\log x}} + O(\frac{c_qx}{q\sqrt{\log x}}(\frac{\log q}{\log x})^{1/7})  \\ \gg \frac{c_qy}{q\sqrt{\log x}} +O(\frac{c_qx}{q\sqrt{\log x}}(\frac{\log\log x}{\log x})^{1/7}) \gg \frac{y}{q\sqrt{\log x}},
\end{multline*}
where in the last inequality we used the fact $y\gg x$ to dismiss the error thanks to it containing the factor $(\log\log x/\log x)^{1/7}$, and therefore, we proved the lemma. 
\section{Proof of The Main Theorem}
We first recall the definition of $T(\alpha, N, s)$ and $R_N(v)$
\[ T(\alpha, N, s) = \frac{1}{N}\sum_{v\in \N}R_N(v)\mathbf{1}_{\{\|v\alpha\|\leq \frac{s}{N}\}} \]
with
\[R_N(v) = \#\{ a_m - a_n = v\mid 1\leq n< m\leq N \}.\]
By the definition of metric Poissonian pair correlation it is enough to prove that for almost every $\alpha$ there are infinitely many $N$'s such that
\[T(\alpha, 10N, 1) \gg_\alpha (\log N)^{1/10}\]
and this is what we aim to obtain, in fact we will be able to prove the above with the constant not depending on $\alpha$ (but the sequence of such $N$'s do).\newline
Let us now consider $a_{10N}$ (the $10N$-th member of the sum of two squares). Landau showed (1908) that the density of the sum of two squares is
\[\frac{N}{a_N} \sim \frac{\kappa}{\sqrt{\log N}}\]
with $0< \kappa < 1$, so for large enough $N$ we have
\[a_{10N}\ge 10N\sqrt{\log N}.\]
Now take any $v\leq 5N\sqrt{\log N}$ for this large enough $N$, then we have
\begin{multline*}
   R_{10N}(v)=\#\{a_m-a_n = v \mid 1\leq n<m\leq 10N\} =\\ \sum_{m\leq a_{10N}-v}b(m)b(m+v) \ge \sum_{m\leq 5N\sqrt{\log N}}b(m)b(m+v). 
\end{multline*}
Now we want to take $n\sim \frac{N}{(\log N)^{1/8}}$ and to examine the sum over $R_{10N}$ along the arithmetic progression $nk$ for $k\leq (\log n)^{5/8}\leq (\log N)^{5/8}$ (so $nk\leq 5N\sqrt{\log N}$ for big enough $N$)
\begin{multline*}
    \sum_{k\leq (\log n)^{5/8}}R_{10N}(nk) \ge \sum_{k\leq (\log n)^{5/8}}\sum_{m\leq 5N\sqrt{\log N}}b(m)b(m+nk)\\ \ge \sum_{k\leq (\log n)^{5/8}}\sum_{m\leq n(\log n)^{5/8}}b(m)b(m+nk) \coloneqq G(n).
\end{multline*}
Now, what we want to do, is to bound $G(n)$ from below for a relatively dense set $B$ of integers $n$, dense enough to have a very good approximations $n\in B$ such that even for $k\leq (\log n)^{5/8}$ we will still have $\|nk\alpha\|\leq \frac{1}{10N}$. 
\newline\newline
We denote by $D_u$ the set of naturals with prime divisors $p\equiv u(4)$, and we set $K = (\log M)^{5/8}$ and $Q = M(\log M)^{5/8}$. Now we use Lemma \ref{lem2} and the fact that $D_1$ is contained within the sum of two squares to obtain
\begin{align*}
    \sum_{n\in [M,2M)}G(n) &= \sum_{M\leq n<2M}\sum_{k\leq (\log n)^{5/8}}\sum_{m\leq n(\log n)^{5/8}}b(m)b(m+nk)\\ &\ge \sum_{m\leq Q}b(m)\sum_{k\leq K}\sum_{M\leq n< 2M}b(m+nk)\\ &\ge \sum_{\substack{m\leq Q\\ m\in D_1}}\sum_{\substack{\frac{1}{2}K\leq k\leq K\\ k\in D_3}}\sum_{\substack{Mk+m\leq l< 2Mk + m\\ l\equiv m(k)}}b(l).
\end{align*}

Since the summation is over $k\in D_1$ and $m\in D_3$ we have that $(m,k) = 1$ and $(m,4) = 1$. Put $x = Mk + m$ and $y = Mk$, since $\frac{K}{2}\leq k \leq K$ and $m\leq Q = MK$ we have that $y \gg x$ and $k\ll \log x \sim \log M$, therefore $m, k, x, y$ satisfy the conditions in Lemma \ref{lem2} so we can use it to estimate the innermost sum of the above inequality's most RHS
\begin{multline*}
\sum_{\substack{Mk+m\leq l< 2Mk + m\\ l\equiv m(k)}}b(l) = \sum_{\substack{x\leq l\leq x+y\\ l\equiv m(k)}}b(l) \\ \gg \frac{y}{k\sqrt{\log x}} = \frac{Mk}{k\sqrt{\log (Mk+m)}} \gg \frac{M}{\sqrt{\log M}}.  
\end{multline*}
Now we use it to see that
\begin{multline*}
    \sum_{n\in[M,2M)}G(n) \gg \sum_{\substack{m\leq Q\\ m\in D_1}}\sum_{\substack{\frac{1}{2}K\leq k\leq K\\ k\in D_3}}\sum_{\substack{Mk+m\leq l< 2Mk + m\\ l\equiv m(k)}}b(l) \\ \gg \frac{M}{\sqrt{\log M}}\sum_{\substack{m\leq Q\\ m\in D_1}}\sum_{\substack{\frac{1}{2}K\leq k\leq K\\ k\in D_3}}1 \gg \frac{M}{\sqrt{\log M}}\frac{K}{\sqrt{\log K}}\frac{Q}{\sqrt{\log Q}}\gg \frac{M^2(\log M)^{2/8}}{\log\log M}.
\end{multline*}
Here we use two facts about the densities of $D_1$ and $D_3$, i.e.
\[ \sum_{\substack{\frac{1}{2}K\leq k\leq K\\ k\in D_3}}1 \gg \frac{K}{\sqrt{\log K}} \]
and
\[\sum_{\substack{m\leq Q \\ m\in D_1}}1\gg \frac{Q}{\sqrt{\log Q}}\]
where in the inequality we took $K = (\log M)^{5/8} $ and $Q = M(\log M)^{5/8}$. These facts can be proven by considering the generating Dirichlet function of $D_1$ and $D_3$, i.e.
\[
\sum_{n\in D_1}\frac{1}{n^\mu}
\]
which is used to derive the asymptotic
\[
\sum_{\substack{n\leq x\\ n\in D_1}}1 \sim \frac{\beta x}{\sqrt{\log x}}
\]
with $\beta = \frac{1}{2\sqrt{2}}> 0$ and this yields the result for
\[\sum_{\substack{\frac{1}{2}K\leq k\leq K\\ k\in D_1}}1\]
by subtracting the sum up to $K/2$ from the sum up to $K$.
The treatment of $D_3$ is similar, and both can be found in \cite{HR} exercise 22 (with $D_1$ containing $2$-divisors but it does not make any difference and the same exact method can be used to obtain the asymptotic for our $D_1$).\newline

This implies that the average of $G$ over $[M,2M)$ is
\[\gg \frac{M(\log M)^{2/8}}{\log\log M}.\]
Now we use the upper bound (see \cite{KF} exercise 2.8)
\[\sum_{m\leq x}b(m)b(m+nk) \ll \prod_{\substack{p|nk\\ p\equiv 3(4)}}(1+\frac{1}{p})\frac{x}{\log x}\]
which holds uniformly for $n$ and $k$ and we apply it with $x = n(\log n)^{5/8}$ to obtain
\begin{multline*}
    G(n) =\sum_{k\leq (\log n)^{5/8}}\sum_{m\leq x}b(m)b(m+nk) \ll \\ \sum_{k\leq (\log n)^{5/8}}\frac{x}{\log x}\prod_{p|nk}(1+\frac{1}{p}) \ll n(\log n)^{2/8}\log\log n
\end{multline*}
with the $\log\log n$ term in the RHS being due to $k \leq \log n$ so using Mertens' inequality and the fact $w(n) \ll \log n$ ($w(n)$ is the number of distinct prime divisors of $n$) we see
\[\prod_{p|nk}(1+\frac{1}{p}) \ll \log\log (nk)\ll \log\log n.\]
The lower bound on the average of $G(n)$ over the interval $[M,2M)$ and the upper bound on all $G(n)$ together imply that
\[
    \#\{n\in [M,2M)\mid G(n)\gg \frac{n(\log n)^{2/8}}{\log\log n}\} \gg \frac{M}{(\log\log M)^2}
\] 
which then implies by Lemma \ref{lem1} that for a.e $\alpha$ the set
\[ B = \{n\mid G(n) \gg \frac{n(\log n)^{2/8}}{\log\log n}\}\] has infinitely many $n$'s such that
\[\|n\alpha\| \leq \frac{1}{n(\log n)^h}\]
for any fixed $h<1$. Now for the final step of the proof, we notice that for any $n$ we can take $N$ with $n\sim \frac{N}{(\log N)^{1/8}}$ so for this choice of $N$ we have for $n\in B$
\[G(n)\gg \frac{n(\log n)^{2/8}}{\log\log n}\gg N(\log N)^{1/10}.\]
Now we want to see that we can find infinitely many $n$'s (or equivalently $N$) such that $\|nk\alpha\| \leq \frac{1}{10N}$ for all $k\leq n(\log n)^{5/8}$ but this is indeed the case because we can apply the above with $h=0.9$ to get
\[\|nk\alpha\|\leq \frac{k}{n(\log n)^h} \leq \frac{1}{n(\log n)^{2/8}} \ll \frac{1}{N(\log N)^{1/8}}\]
and of course for big enough $N$ this implies
\[\|nk\alpha\| \leq \frac{1}{10N}.\] 
Now going back to our definition of $T(\alpha, N, s)$ we get
\[T(\alpha, 10N, 1) = \frac{1}{10N}\sum_{v\in\N}R_{10N}(v)\mathbf{1}_{\{\|v\alpha\|\leq \frac{1}{10N}\}}\]
but for our choice of $n$ and $N$ and due to the computations we did for the sum 
\[\sum_{k\leq (\log n)^{5/8}}R_{10N}(nk)\ge G(n)\]
we have that
\begin{multline*}
    T(\alpha, 10N, 1) = \frac{1}{10N}\sum_{v\in\N}R_{10N}(v)\mathbf{1}_{\{\|v\alpha\|\leq \frac{1}{10N}\}} \ge \\ \frac{1}{10N}\sum_{k\leq (\log n)^{5/8}}R_{10N}(nk) \ge \frac{G(n)}{10N} \gg (\log N)^{1/10} 
\end{multline*}
and now we are done proving there is no pair correlation for almost every $\alpha$.

\end{document}